\documentclass{article}

\usepackage{arxiv}

\usepackage[utf8]{inputenc} % allow utf-8 input
\usepackage[T1]{fontenc}    % use 8-bit T1 fonts
\usepackage{hyperref}       % hyperlinks
\usepackage{url}            % simple URL typesetting
\usepackage{booktabs}       % professional-quality tables
\usepackage{amsfonts}       % blackboard math symbols
\usepackage{nicefrac}       % compact symbols for 1/2, etc.
\usepackage{microtype}      % microtypography
\usepackage{lipsum}		% Can be removed after putting your text content
\usepackage{graphicx}

\usepackage{moreverb}

\usepackage{graphicx, subfigure}

\newcommand\BibTeX{{\rmfamily B\kern-.05em \textsc{i\kern-.025em b}\kern-.08em
T\kern-.1667em\lower.7ex\hbox{E}\kern-.125emX}}

\usepackage{graphicx}
\RequirePackage{fix-cm}
\usepackage{hyperref}
 \usepackage{epsfig} % for postscript graphics files
 \usepackage{amssymb}
 \usepackage{amsmath}
\usepackage{multirow}
\usepackage{epsfig} % for postscript graphics files
\input{epsf}
\usepackage{graphicx}
\usepackage{graphicx, subfigure}
\usepackage[utf8]{inputenc}
\usepackage[english]{babel}

\def\beq{\begin{equation}}
\def\eeq{\end{equation}}
\def\baq{\begin{eqnarray}}
\def\eaq{\end{eqnarray}}
\def\bal{\begin{align} }
\def\eal{\end{align} }
\def\bc{\begin{center}}
\def\ec{\end{center}}
\def\ds{\displaystyle}
\def\fai{\varphi}

\def\fai{\varphi}

\def\gr{\gamma_{r}}
\def\gr1{\gamma_{r1}}

\def\fai{\varphi}

\title{Performance Enhancement of the Recursive Least Squares Algorithms with Rank Two Updates}
\author{ Alexander Stotsky \\
   Systems \& Control \\
Department of Electrical Engineering \\
     Chalmers University of Technology  \\
     Gothenburg  SE - 412 96, Sweden  \\
	\texttt{alexander.stotsky@chalmers.se} \\
        \texttt{alexander.stotsky2@telia.com}  }
\hypersetup{
pdftitle={A template for the arxiv style},
pdfsubject={q-bio.NC, q-bio.QM},
pdfauthor={David S.~Hippocampus, Elias D.~Striatum},
pdfkeywords={First keyword, Second keyword, More},
}
\date{}
\begin{document}
\maketitle

\begin{abstract}
~New recursive least squares algorithms with rank two updates (RLSR2) that
include both exponential and instantaneous forgetting (implemented via a proper choice of the forgetting factor and the window size) are introduced and systematically associated in this report with  well-known RLS algorithms with rank one updates.
Moreover, new properties (which can be used for further performance improvement) of the recursive algorithms associated with the convergence of the inverse of information matrix and parameter vector are established in this report.
The performance  of new algorithms is examined in the problem of estimation of the
grid events in the presence of significant harmonic emissions.
\end{abstract}


\maketitle

\section{Introduction}
\noindent
RLS (Recursive Least Squares) algorithms with forgetting factor are widely used in system identification, signal processing, statistics, control, and in many other applications, \cite{fom}, \cite{lju1}.
% Rls with forgetting
The estimation performance in the weighted least squares problem is highly influenced by the forgetting factor, which discounts exponentially  old measurements and creates a virtual  moving window.
% w + lambda !
\\
The choice of forgetting factor is associated with the trade-off between rapidity and accuracy of estimation.
Introduction of forgetting factor in sliding window, \cite{di}, \cite{card}, \cite{liu} creates extended forgetting mechanism that
includes both exponential and instantaneous forgetting and provides new opportunities for achievement of the trade-off
between rapidity and accuracy.
% rank 1 and 2
\\
The movement of the sliding window is associated with data updating and downdating that results in recursive updates of the information matrix, which can occur sequentially, \cite{c}, \cite{zhang} or simultaneously, \cite{choi}.
Sequential updating and downdating results in computationally complex algorithm, which is difficult to simplify.
Well-known RLS algorithms are associated with updating only and recursive rank one updates, \cite{fom}, \cite{lju1} whereas simultaneous updating/downdating was associated with  computationally efficient
rank two updates in \cite{sto2023}.
\\ This report extends the approach proposed in  \cite{sto2023} and introduces exponential forgetting in the  moving window which prioritizes recent measurements and improves estimation performance for fast varying changes of the signal. The development is performed and associated in a systematic way with well-known RLS algorithms with rank one updates, see Table~1 in \cite{sto2024}.  In addition, rank two gain updates, derived as solution of the least squares problem in
sliding window with exponential forgetting, formed the basis for new Kaczmarz algorithms with improved performance, \cite{kac}.
%RLS algorithms are usually implemented as recursive rank one update of the information matrix, where forgetting
%factor creates a virtual  moving window. The performance of the RLS algorithms with rank two update, \cite{sto2023} which define the window explicitly is improved in this report by introduction of the exponential forgetting factor.
%The method represents a new form of exponential weighting of the data inside of the window which
%prioritizes recent measurements and improves estimation performance for fast varying changes of the signal.
%New RLS algorithms with rank two update which allow to choose both window size and forgetting factor
%are introduced in this report in comparison with RLS algorithms with rank one update, see Table~\ref{ta1}.
%Introduction of two adjustable parameters in new algorithms allows to choose
%sufficiently large window size and reduce forgetting factor for performance improvement in the presence of significant  %measurement noise.  Forgetting factors which are close to one can also be applied for transient performance improvement with %a properly chosen window size when the measurement noise is not significant.
% convergence
\\ Moreover, the gain update $\Gamma_k$ converges to the inverse of the information matrix and the parameters converge to their true values, which is a new property of RLS algorithms with rank two updates discovered in this report, see Section~\ref{glam}.
\\ The response time of the estimation algorithms is restricted by choice of the window size. Short window implies ill-conditioning of the information matrix which results in  sensitivity to numerical calculations and
error accumulation. New convergence properties discovered in  this report allows initialization to approximate inverse
with subsequent convergence of the inverse of information matrix, see Section~\ref{glam}. This property opens new
opportunities for performance improvement in the ill-conditioned case where the difficulties are associated with
matrix inversion. In addition, Newton-Schulz  and Richardson algorithms can be applied for improvements
of the transient performance,\cite{sto2023}, \cite{ifac2022}.
%\\  Moreover, the inverse of the information matrix and the parameter vector can be estimated using
%RLS algorithms with rank two update, which is a new property associated with robustness
%established in this report, see Section~\ref{glam}.
\\ The performance  of new algorithms is examined in the problem of estimation of the
grid events in the presence of significant harmonic emissions  \cite{sto2023},\cite{ifac2022}, \cite{ifac2022b}.
\\ \\
Extended version of this report is presented in \cite{sto2024}:
\\ \\
Stotsky A., Recursive Least Squares Estimation with Rank Two Updates,
Automatika, vol.66, issue 4 , 2025, pp. 619-624 .
\\  \url{https://doi.org/10.1080/00051144.2025.2517431}
\\  \\ see also Kaczmarz projection algorithms with rank two updates and improved performance in \cite{kac}
\\
\url{https://doi.org/10.1007/s11265-024-01915-w}

\section{Least Squares Estimation in Moving Window with Exponential Forgetting}
\label{mw}
\noindent
Estimation of the signal quantities in the moving window is the
most accurate way of monitoring of the wave form distortions
and harmonic emissions in the future electrical networks.
A new form of exponential weighting of the data inside of the window which
prioritizes recent measurements and improves estimation performance for fast varying changes of the  wave form
is considered in this Section.
\subsection*{Problem Formulation and Algorithm Description}
\label{r2}
\noindent
Suppose that a measured oscillating signal can be
presented in the following form ${y}_k = \fai_k^{T} \theta_* $, $k = 1,2,...$ where
the following vector is called the harmonic regressor  $\fai_k^{T} =  [ cos( q_0 k ) ~ sin(q_0 k ) ~ ...~ cos( q_h k ) ~ sin( q_h k ) ]$, where $q_0,...q_h$ are the frequencies  and  $\theta_* $ is
the vector of unknown parameters. The oscillating signal ${y}_k$ is approximated using the model $\hat{y}_k = \fai_k^{T} \theta_k $.
Minimization of the following performance index with exponential forgetting factor $ 0 < \lambda \le 1$ in the moving window
of the size $w$:
\begin{align}
S_k &= \sum_{j= k - (w - 1)}^{k} \lambda^{k-j}~(y_j - \fai^{T}_j \theta_k )^2
\label{ssum}
\end{align}
yields to the following algebraic equations
\begin{align}
 A_k \theta_k &= b_k, ~~  A_k = \sum_{j = k -(w-1)}^{j=k} \lambda^{k-j} \fai_{j}~ \fai_{j}^{T}, ~~  b_k = \sum_{j = k -(w-1)}^{j=k}  \lambda^{k-j} \fai_{j}~ y_j  \label{abk}
 \end{align}
which should be solved with respect to $\theta_k$  in each step $k$.
\\
Notice that the information matrix $A_k$ is defined in (\ref{abk}) as the weighted sum of rank one matrices
can also be defined  as the rank two update
of the matrix  $A_{k-1}$, $k \ge w + 1$. Rank two update is associated with the movement of the window, where
new observation is added (updating) and old observation is deleted (downdating).
In other words, the new data $\fai_{k}$, $y_k$ (with the largest forgetting factor which is equal to one) enter the window and the data with the lowest priority $ \tilde{\fai}_{k-w} = \sqrt{\lambda^w}~\fai_{k-w}$, ~~ $\lambda^w~y_{k-w}$ leave the window in step $k$, \cite{sto2023}, \cite{kac}, \cite{sto2024}:
\begin{align}
 A_k & =  \lambda~A_{k-1} + Q_k~D~Q^{T}_k,~~b_k = \lambda ~ b_{k-1} + d_k \label{r2update}
\end{align}
where $Q_k = [\fai_{k} ~ \tilde{\fai}_{k-w}]$,  $ D = diag[ 1, -1] = \begin{bmatrix}
1 & 0 \\
0 & -1  \end{bmatrix} $ and $d_k =   \fai_{k} ~ y_k -   \lambda^w~\fai_{k-w}~y_{k-w}$.
\\ Notice that the matrix  $Q_k$ contains scaled regressor $ \tilde{\fai}_{k-w} $ in order to avoid singularity
in the case where $\lambda^w \to 0$ for a sufficiently small $\lambda$ and sufficiently large $w$.
Inclusion of sufficiently small $\lambda^w$ in the matrix $D$ (without scaling the regressor)
makes this matrix singular, which results in large estimation errors when calculating the inverse of $A_k$,
\cite{kac}, \cite{sto2024}. The new, minimal/compact form (\ref{r2update}) of the information matrix in the moving window significantly simplifies matrix inversion (in comparison to \cite{zhang}, for example).
\\
Notice also the rank one updates can be obtained as the limiting form of rank two  updates (\ref{r2update}) with
$\lambda^w \to 0$, see Section~\ref{r1}.
\\
The parameter vector in  (\ref{abk}) can be calculated using the
inverse of information matrix, $ \theta_k = A^{-1}_k  b_k$.
Denoting $ \Gamma_{k} = A^{-1}_k $ the recursive update of $ \Gamma_{k} $ via $ \Gamma_{k-1} $
is derived by application of the matrix inversion lemma\footnote{
$ (X + Y W Z)^{-1} = X^{-1} - X^{-1} Y ~ [ W^{-1} + Z~X^{-1}~Y]^{-1} ~ Z X^{-1} $, where $X = \lambda ~A_{k-1} $, $Y = Q_k$, $Z = Q^T_k$ and $W = D$}
to the identity (\ref{r2update}):
\begin{align}
 \Gamma_{k} &= \frac{1}{\lambda} ~ [~ \Gamma_{k-1} - \Gamma_{k-1} ~ Q_k ~ S^{-1} ~ Q^T_k ~ \Gamma_{k-1}~]   \label{gk1}
\end{align} 	
where $ S = \lambda~D + Q^T_k ~ \Gamma_{k-1} ~ Q_k $  is the square {\it capacitance matrix}, \cite{hager}  remains constant for a given window size $w$.
% prog New 2023 Amp_Est_9.m compraisons of different forms !!!
The recursive algorithm for the parameter vector $\theta_k$ is derived  using (\ref{abk}) and  (\ref{gk1}) as follows:
\begin{align}
\theta_k &= A^{-1}_k  ~b_k = \Gamma_{k} ~b_k = \frac{1}{\lambda} ~ [~ \Gamma_{k-1} - \Gamma_{k-1} ~ Q_k ~ S^{-1} ~ Q^T_k ~ \Gamma_{k-1}~] ~
[\lambda b_{k-1} + d_k]  \nonumber \\
    &= \frac{1}{\lambda}  ~ [~ I - \Gamma_{k-1} ~ Q_k ~ S^{-1} ~ Q^T_k ]
~ [\lambda~\underbrace{\Gamma_{k-1}~b_{k-1}}_{\theta_{k-1}} + d_k]  =  [I -   \Gamma_{k-1} ~ Q_k ~ S^{-1} ~ Q^T_k] ~ [ \theta_{k-1} + \Gamma_{k-1} d _k / \lambda ] \label{tet11}
\end{align}
where $I$ is the identity matrix.
The algorithm (\ref{gk1}) and (\ref{tet11}) is initialized as  follows $\Gamma_{w} = A^{-1}_w$ and $A_w~\theta_w = b_w$ and were derived in \cite{sto2023}  for the case $\lambda = 1$.
\\ New algorithm (\ref{gk1}),(\ref{tet11}) provides faster estimation compared to known RLS algorithm (\ref{gks}), (\ref{t1}) for the same forgetting factor.  However, approximately the same fast transient performance can be achieved by reducing the forgetting factor in  known RLS algorithm (\ref{gks}), (\ref{t1}) or the window size in  algorithms described in \cite{sto2023} with forgetting factor which  is equal to one.
\\
Notice that both fast forgetting and small window size imply large condition number of the corresponding information matrix,
sensitivity to measurement noise, numerical operations and significant error accumulation. Algorithm (\ref{gk1}),(\ref{tet11})
has two adjustable parameters (the window size $w$ and the forgetting factor $\lambda$) which provides additional opportunities for optimization (in comparison to (\ref{gks}), (\ref{t1}) and algorithms described in \cite{sto2023})
and hence for performance improvement.
The choice of both forgetting factor and  the window size is associated with the tradeoff between the estimation performance and both sensitivity to measurement noise and  the condition number. On the one hand
the forgetting factor  and the window size should be small enough for fast estimation. On the other hand
the window size should be large enough  and the  forgetting factor should be close to one for small condition number,
reduction of the error accumulation and attenuation of the measurement noise.
Transient performance improvement via reduction of the forgetting factor is preferable in the presence of significant measurement noise. Notice that fast forgetting implies also that  RLS algorithm with rank two updates gets closer to well-known RLS algorithm with rank one updates see Section~\ref{r1}.
\\
Introduction of two adjustable parameters in new algorithm (\ref{gk1}),(\ref{tet11}) allows to choose
sufficiently large window size and reduce forgetting factor for performance improvement in the presence of significant  measurement noise.  Forgetting factors which are close to one can also be applied for transient performance improvement with a properly chosen window size when the measurement noise is not significant.

\subsection*{Rank One Updates as Limiting Form of Rank Two Updates}
\label{r1}
\noindent
Introduction of the forgetting factor allows to establish relationship between RLS
algorithms with  rank two and rank one updates.
Notice that $Q_k$ and $d_k$, see (\ref{r2update}), get the following forms
$Q_k = [\fai_{k} ~ 0]$ and $d_k =   \fai_{k} ~ y_k $, if
$\lambda^w \to 0$ which corresponds to the case of expanding window with the size truncated
by exponential forgetting. Straightforward   substitution of $Q_k$ and $d_k$ in (\ref{gk1}) and (\ref{tet11}) yields to following well-known recursive least squares algorithms:

\begin{align}
 \Gamma_{k} &= \frac{1}{\lambda} ~ [~ \Gamma_{k-1} - \frac{\Gamma_{k-1} ~ \fai_k ~ \fai^T_k ~ \Gamma_{k-1}~}{
 \lambda +  \fai^T_k ~ \Gamma_{k-1} \fai}~]   \label{gks} \\
 \theta_k &= \theta_{k-1} + \frac{\Gamma_{k-1} ~ \fai_k }{
 \lambda +  \fai^T_k ~ \Gamma_{k-1} \fai_k}~~[ y_k - \fai^T_k \theta_k]   \label{t1}
\end{align}

\section{Parameter Calculation with Desired Accuracy \& Simplification of Recursive Matrix Inversion Algorithm}
\label{si}
\noindent
Recursive nature of RLS algorithms (as ideal explicit solution of the system (\ref{abk})
implies error accumulation in finite digit calculations.
Reduction of  the window size and fast forgetting result in ill-conditioned information matrices  and
significant performance deterioration due to error accumulation. Newton-Schulz and Richardson
algorithms can be applied for correction of the inverse of information matrix and parameters, \cite{sto2023}.
The parameter vector in (\ref{abk}) can be calculated with desired accuracy in this case, which essentially improves the estimation performance.
\\
Alternatively, nonrecursive Richardson algorithm described for example in \cite{dub}, \cite{ifac2022} which requires matrix vector multiplications can be used directly  for calculation of the parameters in each step of the moving window:
\begin{align}
\theta_i &= \theta_0 - \sum^{i_*}_{j=0} F_0^j ~~G_0~~ (A_k \theta_0 - b_k), ~~F_0 = I - G_0 A_k
% F_0 = I - \alpha A_k, \alpha = 2 / \| A_k \|_\infty
\label{nonrec1}
\end{align}
via power series expansion until the accuracy requirement is fulfilled. The performance of the algorithm (\ref{nonrec1})
depends on the initial values  $\theta_0$ and $G_0$, where $G_0 = \hat{A}^{-1}_k $ is the estimate of the inverse $A^{-1}_k$ and $\theta_0 = G_0 b_k$. For successful application of this algorithm the approximate inverse $\hat{A}^{-1}_k$
such that $\| F_0 \| <<< 1 $ is required only.
It is clear that the recursive form (\ref{gk1}) with  Newton-Schulz corrections for prevention of the error accumulation can be applied for estimation of $A^{-1}_k$. Notice that the properties of the information matrix in (\ref{abk}) depend on such parameters as forgetting factor and window size. For example, the simplified form (\ref{gks}) with reduced computational complexity associated with rank one updates can be applied for the case where $\lambda^w$ is sufficiently small.
\begin{figure}[t!]
\centerline{\psfig{figure=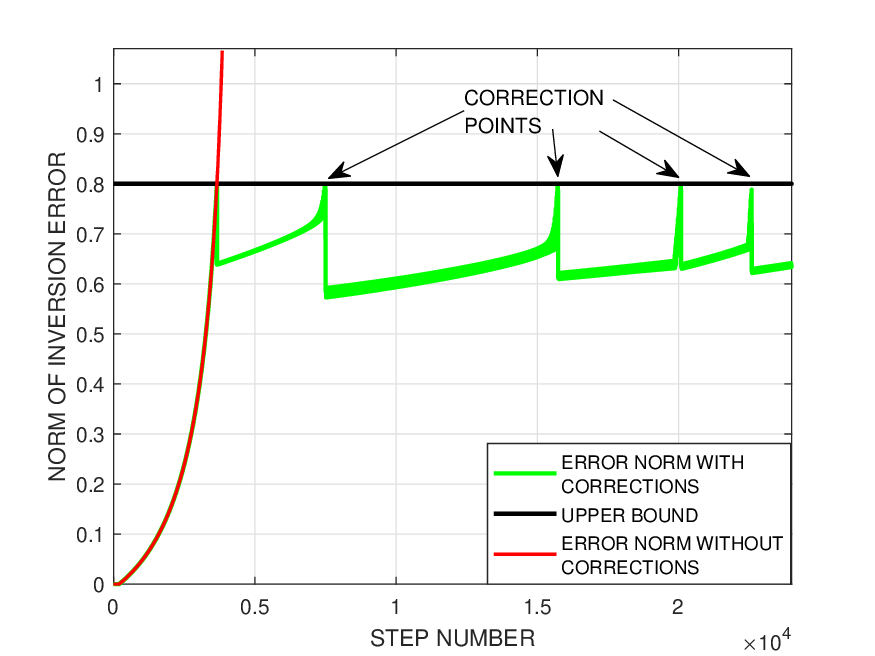,height=75mm}}
\begin{center}
\caption{\small{Matrix inversion error $ \| F_k = I - \Gamma_k A_k \| $ where $\Gamma_k$ is calculated with
 (\ref{gkss}) is plotted with the red line for $w=200$ and $\lambda = 1$. The same error with corrected estimate of the inverse calculated as $\Gamma_{corr~k}  = \Gamma_{k}  + (I - \Gamma_k~ A_k) ~  \Gamma_{k}$ in several points indicated with arrows is plotted with the green line. Pre-specified upper bound of the error norm is plotted with
 the black line.
 }}
\label{ferrac}
\end{center}
\end{figure}
\subsection*{Simplification of Recursive Matrix Inversion Algorithm via Decomposition}
\noindent
Interestingly enough that the form (\ref{gk1}) can be simplified for parallel calculations using the properties of the capacitance matrix $S$, which is the SDD (Strictly Diagonally Dominant) matrix for a sufficiently large window size and sufficiently small forgetting factor. Decomposition of updating and downdating terms in rank two updates can be achieved by neglecting small non diagonal elements of the capacitance matrix. Explicit evaluation of the diagonal elements of the capacitance matrix with subsequent substitution in (\ref{gk1}) results in the following recursive equation (decomposed on updating and downdating terms for parallel calculations) for approximation of inverse of information matrix $ \Gamma_{k} \approx A^{-1}_k $ :
\begin{align}
 \Gamma_{k} &= \frac{1}{\lambda} ~ \{  \Gamma_{k-1}~ - \underbrace{[~\frac{\Gamma_{k-1} ~ \fai_k ~ \fai^T_k ~ \Gamma_{k-1}~}{\lambda +  \fai^T_k ~ \Gamma_{k-1} \fai_k}~]}_{\text{updating}}
 -   \underbrace{[~ \frac{\Gamma_{k-1} ~ \tilde{\fai}_{k-w} ~ \tilde{\fai}^{T}_{k-w} ~ \Gamma_{k-1}~}{
 - \lambda +  \tilde{\fai}^{T}_{k-w} ~ \Gamma_{k-1} \tilde{\fai}^{T}_{k-w}} ~ ]}_{\text{downdating}} \}     \label{gkss}
\end{align}
Notice that the denominators in updating and downdating terms are constant and (\ref{gkss}) has the limiting form of rank one updates  (\ref{gks}) when  $\lambda^w \to 0$ and downdating
term disappears.
% program for error accumulation Amp_Est_8.m !!!
\\ The drawback of this approximate relation is error accumulation and significant deterioration of the inversion accuracy for a sufficiently large window sizes and forgetting factors which are close to one. The error accumulation  problem is illustrated in Figure~\ref{ferrac}, where the infinity matrix norm of the inversion error  $ \| F_k = I - \Gamma_k A_k \| $ is plotted with the red line. The convergence rate of the Richardson algorithm (\ref{nonrec1}) strongly depends on accuracy of the estimate of the inverse $A^{-1}_k$ measured by this infinity norm. The error can be corrected using the following one step Newton–Schulz algorithm $\Gamma_{corr~k}   = \Gamma_{k}  + (I - \Gamma_k ~A_k)~ \Gamma_{k}$ which requires two matrix products only and can be implemented using parallel calculations.
Figure~\ref{ferrac} shows that the  norm of the inversion error does not exceed the pre-specified value $0.8$ for system with corrections that guarantees fast convergence of the Richardson algorithm (\ref{nonrec1}).
Notice that the error accumulation problem for the algorithm (\ref{gkss}) can be eliminated/reduced by reduction of the forgetting factor.

\section{Parameter Estimation with RLS Algorithm with Rank Two Updates}
\label{glam}
\noindent
\subsection{Description of Algorithms and Error Models}
\label{glam1}
\noindent
\begin{figure}[t!]
\centerline{\psfig{figure=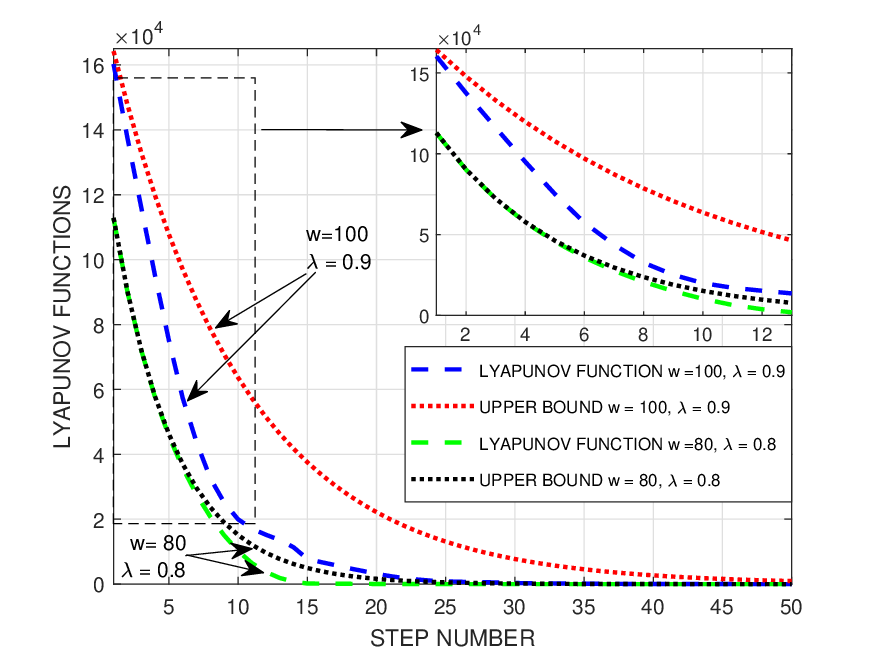,height=75mm}}
\begin{center}
\caption{\small Lyapunov functions  $V_k$  (green and blue dashed lines) are plotted with the upper bounds
 $V_k \le \lambda^{k} ~ V_0 $ (black and red dotted lines) for $w=80$, $\lambda = 0.8$ and $w=100$, $\lambda = 0.9$
 respectively. }
\label{lyap11}
\end{center}
\end{figure}

The algorithms (\ref{gk1}), (\ref{tet11}) can be written in the following form, \cite{kac}, \cite{sto2024} :
\begin{align}
 \Gamma_{k} &= \frac{1}{\lambda} ~ [~ \Gamma_{k-1} - \Gamma_{k-1} ~ Q_k ~ S^{-1} ~ Q^T_k ~ \Gamma_{k-1}~]   \label{gk2} \\
  \theta_k  &= \theta_{k-1} - \Gamma_{k-1}  ~ Q_k ~ S^{-1} ~ [ Q^T_k~ \theta_{k-1} - \tilde{y}_k ] \label{ts1}
\end{align}
where $\tilde{y}_k $ is the augmented output, provided that the matrix
$Q^T_k ~ \Gamma_{k-1} ~ Q_k $ is invertible.
\\ For system  (\ref{r2update}), (\ref{gk2}) and (\ref{ts1}) the following error model is valid,\cite{sto2024} :

\begin{align}
   E_k &= ( I - \Gamma_{k-1}  ~ Q_k ~ S^{-1} ~ Q^T_k)~ E_{k-1}    \label{erg1} \\
 \tilde{\theta}_k &= ( I - \Gamma_{k-1}  ~ Q_k ~ S^{-1} ~ Q^T_k ) ~ \tilde{\theta}_{k-1} \label{tilte1}
\end{align}
where $E_k = I - \Gamma_{k}~ A_k$ and $\tilde{\theta}_k = \theta_k  - \theta_* $ are matrix inversion and parameter estimation errors.

\subsection{Convergence Properties of Matrix Inversion and Parameter Errors}
\label{glam2}
\noindent
The convergence of the matrix inversion error can be established  by explicit evaluation of $E_k$ along the solution of (\ref{erg1}). Notice that similar convergence of the matrix inversion error is valid for RLS algorithm with rank one updates.
\\
Transient parameter estimation performance is evaluated (using arguments similar to \cite{st2012}, \cite{st2013})
by considering the first difference
$V_k - V_{k-1}$ of the Lyapunov function $V_k =  \tilde{\theta}_k^T ~ A_k ~ \tilde{\theta}_k$  under the assumption that $E_k = 0$:

\begin{align}
V_k - V_{k-1} & = - \lambda~\tilde{\theta}^T_{k-1}~ Q_k ~ S^{-1} ~ Q^T_k~\tilde{\theta}_{k-1} - (1 - \lambda)~V_{k-1}   \label{vk1}
\end{align}
which implies that $\ds V_k \le \lambda^{k} ~ V_0 $ and $\ds  \| \tilde{\theta}_k \| \le \sqrt{\frac{\lambda^{k} ~ V_0}{\lambda_{min}~(A_k) }} $ provided that $\tilde{\theta}^T_{k-1}~ Q_k ~ S^{-1} ~ Q^T_k~\tilde{\theta}_{k-1} \ge 0 $.
Lyapunov functions and their upper bounds for different window sizes and forgetting factors are presented in
Figure~\ref{lyap11}, which shows that Lyapunov approach provides relatively tight bounds on estimation errors.  The convergence analysis without the assumption that $E_k = 0$ is presented in \cite{sto2024}.

\section{Conclusion}
\label{conc}
\noindent
The forgetting factor which allows  exponential weighting of the data inside of the moving window,
prioritizes recent measurements and improves estimation performance for fast varying changes of the signal
was introduced in RLS algorithms with rank two updates.
It is shown on the problem of the estimation of the grid events  that a proper choice of two adjustable parameters (window size and forgetting factor) in new algorithms essentially improves estimation performance.
New RLS algorithms with rank two updates were systematically associated with well-known RLS algorithms with rank one updates,
see Table~1 in \cite{sto2024}.
Finally, new properties (which can be used for further performance improvement) of the recursive algorithms associated with the convergence of the inverse of information matrix and parameter vector were established in this report.

%\begin{ack}
%Place acknowledgments here.
%\end{ack}

%\bibliography{ifacconf}             % bib file to produce the bibliography
                                                     % with bibtex (preferred)

\newpage

%\section{Appendix}
%\label{appen}

\end{document}